\documentclass[reqno,12pt]{amsart}
\theoremstyle{plain}    
\newtheorem{thm}{Theorem} 
\theoremstyle{plain}    
\newtheorem{prop}[thm]{Proposition}

\newtheorem{remark}[thm]{Remark}

\usepackage{amscd,amssymb,comment,epic,eufrak,euscript,graphics}

\newcount\theTime
\newcount\theHour
\newcount\theMinute
\newcount\theMinuteTens
\newcount\theScratch
\theTime=\number\time
\theHour=\theTime
\divide\theHour by 60
\theScratch=\theHour
\multiply\theScratch by 60
\theMinute=\theTime
\advance\theMinute by -\theScratch
\theMinuteTens=\theMinute
\divide\theMinuteTens by 10
\theScratch=\theMinuteTens
\multiply\theScratch by 10
\advance\theMinute by -\theScratch

\def\today{{\number\day\space
 \ifcase\month\or
  January\or February\or March\or April\or May\or June\or
  July\or August\or September\or October\or November\or December\fi
 \space\number\year}}

\newcommand\At{{\widetilde A}}

\newcommand\Cpx{{\mathbf C}}

\newcommand\diag{\text{\rm diag}}

\newcommand\Ht{{\widetilde H}}

\newcommand\id{{\operatorname{id}}}

\newcommand\Ints{{\mathbf Z}}

\newcommand\FEu{{\EuScript F}}

\newcommand\KEu{{\EuScript K}}                   

\newcommand\LEu{{\EuScript L}}                   

\newcommand\lspan{\mathrm{span}\,}

\newcommand\Nats{{\mathbf N}}

\newcommand\Tt{{\widetilde T}}

\newcommand{\vp}{\varepsilon}

\begin{document}

\pagestyle{myheadings}

 \title{The completely bounded approximation property for extended 
Cuntz--Pimsner 
algebras}

 \author{Kenneth J.\ Dykema and Roger R.\ Smith}

 \address{\hskip-\parindent
 Department of Mathematics \\
 Texas A\&M University \\
 College Station TX 77843--3368, USA}
 \email{kdykema@math.tamu.edu}
 \email{rsmith@math.tamu.edu}

 \thanks{The first named author was supported in part by a grant from the NSF}

 \date{11 November, 2003}

 \begin{abstract}
The extended Cuntz--Pimsner algebra $E(H)$, introduced by Pimsner, is 
constructed from a Hilbert 
$B,B$--bimodule $H$ over a $C^*$--algebra $B$. In this paper we investigate 
the Haagerup invariant $\Lambda(\cdot)$ for these algebras, the main result 
being that 
\[\Lambda(E(H))=\Lambda(B)\]
when $H$ is  full over $B$. In particular, $E(H)$ has the 
completely bounded approximation property if and only if the same is true for 
$B$.
 \end{abstract}

 \maketitle

 \markboth{\sc K.J.\ Dykema, R.R.\ Smith}{\sc The completely bounded approximation property}

\section{Introduction}

The Haagerup invariant for a $C^*$--algebra $A$ is defined to be the smallest 
constant $\Lambda(A)$ for which there exists a net $\{\phi_\alpha\colon A\to 
A\}_{\alpha \in N}$ of finite rank maps satisfying
\begin{equation}\label{eq1.1}
\lim_\alpha \,\|\phi_\alpha(x)-x\|=0,\ x\in A,
\ \ {\mathrm{ and }}\ \ \|\phi_\alpha\|_{cb}\leq \Lambda(A),\ \alpha \in N.
\end{equation}
If no such net exists then $\Lambda(A)$ is defined to be $\infty$, while if 
$\Lambda(A)<\infty$ then $A$ is said to have the {\it{completely bounded 
approximation property}} (CBAP). The definition of $\Lambda(A)$ arose from 
\cite{H,deCan}. In the first of these it was shown that $C^*_r(\mathbb F_2)$, 
the reduced $C^*$--algebra of the free group on two generators, has such a net 
of contractions, and a stronger result using complete contractions was 
obtained in the second paper. Subsequently, many examples of different values 
of $\Lambda(\cdot)$ were given in \cite{CowH}. An interesting problem is to 
investigate the behavior of $\Lambda(\cdot)$ under the standard constructions 
of $C^*$--algebra theory. In \cite{SScam}, the formula $\Lambda(A_1\otimes 
A_2)=\Lambda(A_1)\Lambda(A_2)$ was established for the minimal tensor product, 
while $\Lambda(A\rtimes_\alpha G)=\Lambda(A)$ was proved for discrete amenable 
groups in \cite{SS} and for general amenable groups in \cite{NS}. Our 
objective in this paper is to show that $\Lambda(B)=\Lambda(E(H))$, where 
$E(H)$ is the extended Cuntz--Pimsner algebra arising from a 
 Hilbert $B,B$--bimodule over a $C^*$--algebra $B$, \cite{P}. The 
$C^*$--algebras $E(H)$ appear in several areas of operator algebras, notably
in the work of Muhly and Solel, \cite{MS1,MS2}, and in work of the first 
author with Shlyakhtenko, \cite{DS}. It was shown in the latter paper that 
$E(H)$ is exact if and only if $B$ is exact. Since exactness is a well known 
consequence of the CBAP, (see, for example, the argument in \cite[Theorem 
3.1(vii)]{NS}), this suggested the connection between $\Lambda(E(H))$ 
and $\Lambda(B)$.
(See also the remarks at the end of this paper.)
We now review briefly the definition of $E(H)$.

A right Hilbert $B$--module $H$ has a $B$--valued inner product $\langle\cdot 
,\cdot\rangle_B$, conjugate linear and linear respectively in the first and 
second variables, and is said to be {\it{full}} if $\{\langle h_1 ,h_2 
\rangle_B \colon h_1,\,h_2\in H\}$ generates $B$. The $C^*$--algebra $\LEu(H)$ 
consists of the $B$--linear operators $T\colon H\to H$ for which  
there is a $B$--linear $T^*\colon H\to H$ satisfying $\langle Th_1 ,h_2 
\rangle_B=\langle h_1 
,T^*h_2 \rangle_B$, and operators in $\LEu(H)$ are called {\it{adjointable}}.
The $C^*$--algebra $\LEu(H)$ contains the closed ideal $\KEu(H)$, that is 
generated by the 
maps of the form
\[\theta_{x,y}(h)=x\langle y, h\rangle_B,\ \ \ h\in H,\]
for arbitrary pairs $x,y\in H$, \cite{L}. If there is an injective 
$*$--homomorphism $\rho 
\colon B\to \LEu(H)$, then there is a left action of $B$ on $H$ by 
$(b,h)\mapsto \rho(b)h$, and we say that $H$ is a $B,B$--bimodule. The full 
Fock space $\FEu(H)$ is defined to be $B\oplus\bigoplus_{n\geq 1}H^{(\otimes 
B)n}$, where $H^{(\otimes B)n}$ is the $n$--fold tensor product 
$H\otimes_BH\otimes_B\cdots\otimes_BH$, which is also a Hilbert 
$B,B$--bimodule. For $h\in H$, the operator $\ell(h)\colon \FEu(H)\to \FEu(H)$ 
is defined on generators by 
\begin{align*} 
\ell(h)h_1\otimes\cdots\otimes h_n&=h\otimes h_1\otimes\cdots\otimes h_n,
\ \ \ h_i \in H,\\
\ell(h)b&=hb.\ \ \ b\in B.
\end{align*}
These are bounded adjointable operators on $\FEu (H)$ and satisfy
\begin{align*}
\ell(h_1)^*\ell(h_2)&=\langle h_1 ,h_2 \rangle_B,\ \ \ h_i \in H,\\
b_1\ell(h)b_2&=\ell(b_1hb_2),\ \ \ b_i \in B.
\end{align*}
Then $E(H)$ is the $C^*$--algebra generated by $\{\ell(h)\colon h\in H\}$, 
introduced by Pimsner in \cite{P}.

Our approach to investigating $\Lambda(E(H))$ follows the methods of 
\cite{DS}. There, a sequence of operations was given to construct $E(H)$ from 
$B$ in such a way that exactness was preserved at each step. Here we show that 
the Haagerup invariant is preserved, with the main technical device being 
Theorem \ref{thm:Lses} concerning quotients of $C^*$--algebras.

\section{Results}

For a short exact sequence
\[0\to J \to A \to A/J \to 0\]
of $C^*$--algebras, we cannot expect, in general, a relationship between 
$\Lambda(A)$, $\Lambda(J)$ and $\Lambda(A/J)$, \cite[Section 5]{NS}. However, 
our first result shows 
that these quantities are closely linked when the short exact sequence splits.
The proof requires the notion of a {\it{quasi--central approximate identity}} 
for a closed ideal $J$ contained in a $C^*$--algebra $A$, introduced in 
\cite{AP,A}. This is an approximate identity $\{e_\alpha\}_{\alpha \in N}$,
$0\leq e_\alpha \leq 1$, which has the additional property of asymptotically 
commuting with the elements of $A$ in the sense that
\[\lim_{\alpha}\,\|e_\alpha x-xe_\alpha\|=0,\ \ \ x\in A.\]
Such approximate identities always exist, \cite{AP,A}.

\begin{thm}\label{thm:Lses}
Let $J$ be an ideal in a $C^*$--algebra $A$, let $\pi\colon \ A\to A/J$ be the 
quotient map, and suppose that there exists a completely contractive map 
$\rho\colon \ A/J\to A$ such that $\pi \rho = \id_{A/J}$. Then
\[
\Lambda(A) = \max\{\Lambda(J), \Lambda(A/J)\}.
\]
\end{thm}

\begin{proof}
Fix $a_1,\ldots, a_n\in A$, $\|a_i\|\le 1$, and fix $\vp>0$. We will construct 
a 
finite rank map $\gamma\colon \ A\to A$ such that $\|\gamma(a_i) - a_i\| < 
\vp$, 
$1\le i\le n$, and $\|\gamma\|_{cb} \le \max\{\Lambda(J), \Lambda(A/J)\}$. 
This 
will then prove ``$\le$'' in the equality.

Fix $\delta>0$, to be chosen later. Let $\{e_\alpha\}_{\alpha\in N}$ be a 
quasi-central approximate identity for $J$ satisfying $0\le e_\alpha \le 1$. 
By 
definition of $\Lambda(A/J)$, there exists a finite rank map $\psi\colon \ 
A/J\to A/J$, $\|\psi\|_{cb}\le \Lambda(A/J)$, such that
\begin{equation*}
\|\psi(\pi(a_i)) - \pi(a_i)\| < \delta
\end{equation*}
for $1\le i\le n$. Define a finite rank map $\widetilde\psi\colon \ A\to A$ by
\begin{equation*}
\widetilde\psi(x) = \rho(\psi(\pi(x))),\qquad x\in A.
\end{equation*}
Clearly $\|\widetilde\psi\|_{cb} \le \Lambda(A/J)$ and, since 
$\pi(\widetilde\psi(a_i)) = \psi(\pi(a_i))$, we may choose elements $j_i\in J$ 
such that
\begin{equation}\label{eq3}
\|\widetilde\psi(a_i) - (a_i+j_i)\| < \delta,\qquad 1\le i\le n.
\end{equation}
Now $\{e_\alpha\}_{\alpha\in N}$ is an approximate identity for $J$ so, for 
each 
$\alpha\in N$, we may choose $\beta(\alpha)\in N$, $\beta(\alpha)>\alpha$, 
such 
that
\begin{equation}\label{eq4}
\|e_{\beta(\alpha)} e^{1/2}_\alpha -e ^{1/2}_\alpha\|<\delta.
\end{equation}
Note that taking adjoints in \eqref{eq4} gives 
\begin{equation}\label{eq5}
\|e^{1/2}_\alpha e_{\beta(\alpha)} - e^{1/2}_\alpha\| < \delta
\end{equation}
also. For each $\alpha\in N$, choose a finite rank map $\phi_\alpha\colon \ 
J\to 
J$, $\|\phi_\alpha\|_{cb}\le \Lambda(J)$ such that
\begin{equation}\label{eq6}
\|\phi_\alpha(e_{\beta(\alpha)}a_ie_{\beta(\alpha)}) - e_{\beta(\alpha)} a_i 
e_{\beta(\alpha)} \| < \delta,\qquad 1\le i\le n.
\end{equation}
Then define $\widetilde\phi_\alpha\colon \ A\to A$ by
\begin{equation*}
\widetilde\phi_\alpha(x) = \phi_\alpha(e_{\beta(\alpha)}xe_{\beta(\alpha)}), 
\qquad x\in A.
\end{equation*}
These maps are finite rank, and satisfy $\|\widetilde\phi_\alpha\|_{cb}\le 
\Lambda(J)$. Then, for each $\alpha\in N$, define $\gamma_\alpha\colon \ A\to 
A$ 
by
\begin{equation}\label{eq8}
\gamma_\alpha(x) = e^{1/2}_\alpha \widetilde\phi_\alpha(x) e^{1/2}_\alpha + 
(1-e_\alpha)^{1/2} \widetilde\psi(x)(1-e_\alpha)^{1/2},\qquad x\in A.
\end{equation}
Each of these maps can be expressed as a matrix product
\begin{equation*}
\gamma_\alpha(x) = (e^{1/2}_\alpha,
 (1-e_\alpha)^{1/2}) \left(\begin{matrix} \widetilde\phi_\alpha(x)&0\\ 
0&\widetilde\psi(x)\end{matrix}\right) \left(\begin{matrix} e^{1/2}_\alpha\\ 
(1-e_\alpha)^{1/2}\end{matrix}\right)
\end{equation*}
from which it is clear that $\|\gamma_\alpha\|_{cb} \le \max\{\Lambda(J), 
\Lambda(A/J)\}$. It remains to be shown that
\begin{equation}\label{eq10}
\|\gamma_\alpha(a_i) - a_i\|<\vp, \qquad 1\le i\le n,
\end{equation}
for a sufficiently large choice of $\alpha\in N$. We will estimate the two 
terms 
on the right hand side of \eqref{eq8} separately, when $x$ is replaced by 
$a_i$.

The term $e^{1/2}_\alpha\widetilde\phi_\alpha(a_i) e^{1/2}_\alpha$ equals, by 
definition, $e^{1/2}_\alpha\phi_\alpha(e_{\beta(\alpha)}a_i 
e_{\beta(\alpha)}) e^{1/2}_\alpha$ so, from \eqref{eq6},
\begin{equation*}
\|e^{1/2}_\alpha\widetilde\phi_\alpha(a_i)e^{1/2}_\alpha - e^{1/2}_\alpha 
e_{\beta(\alpha)} a_i e_{\beta(\alpha)} e^{1/2}_\alpha\| < \delta
\end{equation*}
for $1\le i\le n$. A simple triangle inequality argument, using \eqref{eq4} 
and 
\eqref{eq5}, then shows that
\begin{equation}\label{eq12}
\|e^{1/2}_\alpha \widetilde\phi_\alpha(a_i) e^{1/2}_\alpha - e^{1/2}_\alpha 
a_ie^{1/2}_\alpha\| < 3\delta,
\end{equation}
for $1\le i\le n$.

The second term, $(1-e_\alpha)^{1/2} \widetilde\psi(a_i) (1-e_\alpha)^{1/2}$, 
in 
\eqref{eq8} is within $\delta$ of $(1-e_\alpha)^{1/2}(a_i + j_i) 
(1-e_\alpha)^{1/2}$, from \eqref{eq3}. Combining this estimate with 
\eqref{eq8} 
and \eqref{eq12} gives 
\begin{equation*}
\|\gamma_\alpha(a_i) - e^{1/2}_\alpha a_ie^{1/2}_\alpha - (1-e_\alpha)^{1/2} 
a_i(1-e_\alpha)^{1/2}\| \le \|(1-e_\alpha)^{1/2} j_i(1-e_\alpha)^{1/2}\| + 
3\delta, \quad 1\le i\le n.
\end{equation*}
Thus
\begin{align}
\|\gamma_\alpha(a_i) -a_i\| &\le \|a_i-e^{1/2}_\alpha a_ie^{1/2}_\alpha - 
(1-e_\alpha)^{1/2} a_i(1-e_\alpha)^{1/2}\|\nonumber\\
\label{eq14}
&\quad + \|(1-e_\alpha)^{1/2} j_i(1-e_\alpha)^{1/2}\| + 3\delta, \qquad 1\le 
i\le n.
\end{align}
The nets $\{e^{1/2}_\alpha\}_{\alpha\in N}$, and 
$\{(1-e_\alpha)^{1/2}\}_{\alpha\in N}$ asymptotically commute with each $a_i$, 
and
\begin{equation*}
\|(1-e_\alpha)^{1/2} j_i\|^2 = \|j^*_i(1-e_\alpha)j_i\|\to 0
\end{equation*}
as $\alpha\to\infty$, so a sufficiently large choice of $\alpha$ in 
\eqref{eq14} gives
\begin{equation*}
\|\gamma_\alpha(a_i)-a_i\| < 4\delta,\qquad 1\le i\le n.
\end{equation*}
Now take $\delta$ to be $\vp/4$, and \eqref{eq10} is proved.

We now turn to the inequality $\Lambda(J)\le\Lambda(A)$. Given $j_1,\ldots, 
j_n\in J$, $\|j_i\|\le 1$, and $\vp>0$, we may choose a finite rank completely 
bounded map $\phi\colon \ A\to A$ such that $\|\phi\|_{cb}\le \Lambda(A)$ and
\begin{equation*}
\|\phi(j_i)-j_i\| <\vp,\qquad 1\le i\le n.
\end{equation*}
For each $\alpha\in N$, define $\phi_\alpha\colon \ J\to J$ by
\begin{equation*}
\phi_\alpha(j) = e_\alpha\phi(j)e_\alpha,\qquad 1\le i\le n.
\end{equation*}
then each $\phi_\alpha$ has finite rank, $\|\phi_\alpha\|_{cb} \le 
\Lambda(A)$, 
and, for $1\le i\le n$,
\begin{equation}\label{eq19}
\|\phi_\alpha(j_i)-j_i\| < \|e_\alpha j_ie_\alpha - j_i\| + \vp.
\end{equation}
A sufficiently large choice of $\alpha$ in \eqref{eq19} gives
\begin{equation*}
\|\phi_\alpha(j_i) - j_i\| <\vp,\qquad 1\le i\le n,
\end{equation*}
showing that $\Lambda(J) \le \Lambda(A)$.

Finally we show that $\Lambda(A/J) \le \Lambda(A)$. Consider elements 
$\pi(a_1),\ldots, \pi(a_n)\in A/J$, and fix $\vp>0$. There exists $\phi\colon 
\ 
A\to A$ such that $\phi$ has finite rank, $\|\phi\|_{cb}\le \Lambda(A)$, and 
\begin{equation}\label{eq21}
\|\phi(\rho\pi(a_i)) - \rho\pi(a_i)\| <\vp
\end{equation}
for $1\le i\le n$. Apply $\pi$ to \eqref{eq21} to obtain
\begin{equation*}
\|\pi\phi\rho(\pi(a_i)) - \pi(a_i)\| <\vp
\end{equation*}
for $1\le i\le n$. Now define $\widetilde\phi\colon \ A/J\to A/J$ by
\begin{equation*}
\widetilde\phi = \pi\phi\rho.
\end{equation*}
Clearly $\widetilde\phi$ has finite rank, $\|\widetilde\phi\|_{cb}\le 
\Lambda(A)$, and
\begin{equation*}
\|\widetilde\phi(\pi(a_i)) - \pi(a_i)\| <\vp
\end{equation*}
for $1\le i\le n$. This shows that $\Lambda(A/J) \le \Lambda(A)$, proving the 
result.
\end{proof}

\begin{remark}\rm
If we loosen the hypotheses of Theorem~\ref{thm:Lses} and require only that 
$\rho$ be completely bounded,
then the same proof yields the inequalities
\[
\Lambda(A)\le\max\{\Lambda(J),\|\rho\|_{cb}\,\Lambda(A/J)\},\qquad
\Lambda(J)\le\Lambda(A),\qquad
\Lambda(A/J)\le\|\rho\|_{cb}\,\Lambda(A).
\]
\end{remark}

The next three results are preparatory for Theorem \ref{main}, and will handle 
some technical points arising there.  

\begin{prop}\label{prop:ST}
Let $A$ and $B$ be $C^*$--algebras.
Let $E$ be a right Hilbert $A$--module and let $F$ be a right Hilbert 
$B$--module
with a $*$--homomorphism $\phi:A\to\LEu(F)$.
Consider the internal tensor product $E\otimes_\phi F$, which is a right 
Hilbert $B$--module.
Let $S\in\LEu(E)$ and suppose $T\in\LEu(F)$ is such that $T$ and $\phi(a)$ 
commute for all $a\in A$.
Then there is $R\in\LEu(E\otimes_\phi F)$ satisfying $R(e\otimes 
f)=(Se)\otimes(Tf)$;
we will write $R=S\otimes T$.
\end{prop}
\begin{proof}
The operator $S\otimes\id_F$ is well known to belong to $\LEu(E\otimes_\phi 
F)$;
(see~\cite[p.\ 42]{L}).
It will suffice to show $\id_E\otimes T\in\LEu(E\otimes_\phi F)$, for in 
general we will have
$S\otimes T=(S\otimes\id_F)\circ(\id_E\otimes T)$.
Hence, without loss of generality, we assume $S=\id_E$.

The $\Cpx$--linear map from the algebraic tensor product (over $\Cpx$) 
$E\otimes F$ to itself
defined by $e\otimes f\mapsto e\otimes Tf$ is a right $B$--module map and 
sends 
the
submodule 
\[
N=\lspan\{ea\otimes f-e\otimes\phi(a)f\mid e\in E,\,f\in F,\,a\in A\}
\]
into itself.
In order to see that the resulting map $(E\otimes F)/N\to(E\otimes F)/N$ gives 
rise to a bounded
$B$--linear map $E\otimes_\phi F\to E\otimes_\phi F$, let $n\in\Nats$,
$e_1,\ldots,e_n\in E$ and $f_1,\ldots, f_n\in F$, and let $g=\sum_{i=1}^n 
e_i\otimes f_i$.
Then
\[
\langle g,g\rangle_B=\sum_{i,j}\langle f_i,\phi(\langle 
e_i,e_j\rangle_A)f_j\rangle_B
=\langle 
f,\phi_n(X)f\rangle_{F^n}=\langle\phi_n(X^{1/2})f,\phi_n(X^{1/2})f\rangle_{F^n},
\]
where $f$ is the column vector $(f_1,\ldots,f_n)^t$ in the Hilbert $B$--module 
$F^n$,
where $X$ is the matrix $X=(\langle e_i,e_j\rangle_A)_{1\le i,j\le n}\in 
{\mathbb M}_n(A)$,
which, by~\cite[Lemma 4.2]{L}, is positive, and where $\phi_n:{\mathbb 
M}_n(A)\to 
{\mathbb M}_n(\LEu(F))=\LEu(F^n)$
is the $*$--homomorphism obtained by application of $\phi$ to each element of 
a matrix.
On the other hand, letting $h=\sum_{i=1}^n e_i\otimes Tf_i$,
we have
\begin{align*}
\langle h,h\rangle_B&=\sum_{i,j}\langle Tf_i,\phi(\langle 
e_i,e_j\rangle_A)Tf_j\rangle_B \\
&=\langle\phi_n(X^{1/2})\Tt f,\phi_n(X^{1/2})\Tt f\rangle_{F^n}
=\langle\Tt\phi_n(X^{1/2})f,\Tt\phi_n(X^{1/2})f\rangle_{F^n},
\end{align*}
where $\Tt=\diag(T,\ldots,T)\in {\mathbb M}_n(\LEu(F))$.
Therefore, ({\em cf}\/~\cite[Prop.\ 1.2]{L}),
\[
\langle 
h,h\rangle_B=|\Tt\phi_n(X^{1/2})f|^2\le\|\Tt\|^2|\phi_n(X^{1/2})f|^2=\|T\|^2\langle 
g,g\rangle_B
\]
and $\|h\|\le\|T\|\,\|g\|$.
Consequently, we get $\id_E\otimes T:E\otimes_\phi F\to E\otimes_\phi F$ with 
$\|\id_E\otimes T\|\le\|T\|$.
An easy calculation shows that $\id_E\otimes T^*$ is the adjoint of 
$\id_E\otimes T$.
Thus $\id_E\otimes T\in\LEu(E\otimes_\phi F)$.
\end{proof}

\begin{prop}\label{prop:Lind}
Let $A_1\subseteq A_2\subseteq\cdots\subseteq A$
be an increasing chain of $C^*$--subalgebras of a $C^*$--algebra $A$, such 
that $\bigcup_{n=1}^\infty A_n$
is dense in $A$.
Suppose that there are conditional expectations $\rho_n:A_{n+1}\to A_n$ onto 
$A_n$, $n\geq 1$.
Then
\begin{equation}\label{eq:LAn}
\Lambda(A)=\sup_{n\in\Nats}\Lambda(A_n).
\end{equation}
\end{prop}
\begin{proof}
From the family $(\rho_n)_{n\ge1}$ we obtain conditional expectations 
$\psi_n:A\to A_n$ onto $A_n$, $n\geq 1$. It then follows that  
$\Lambda(A)\ge\Lambda(A_n)$ for all $n\geq 1$.

To see the reverse inequality in~\eqref{eq:LAn}, let $n\in\Nats$, let 
$F\subseteq A_n$
be a finite subset and  let $\varepsilon>0$.
Then there is a finite rank map $\phi:A_n\to A_n$ such that 
$\|\phi(x)-x\|<\varepsilon$ for all
$x\in F$ and $\|\phi\|_{cb}\le\Lambda(A_n)$.
But $\phi\circ\psi_n:A\to A_n\subseteq A$ satisfies 
$\|\phi\circ\psi_n(x)-x\|<\varepsilon$ for all $x\in F$
and $\|\phi\circ\psi_n\|_{cb}=\|\phi\|_{cb}$.
\end{proof}

\begin{prop}\label{prop:KH}
Let $B$ be a $C^*$--algebra and let $H$ be a 
right Hilbert
$B$--module such that $\{\langle h_1,h_2\rangle_B\colon h_1,h_2\in H\}$ 
generates 
$B$.
Then $\Lambda(\KEu(H))\le\Lambda(B)$.
\end{prop}
\begin{proof}
Following the notation of \cite{B}, let $C_n(B)$ denote the right Hilbert 
$B$--module which consists of columns over $B$ of length $n$ with the 
$B$--valued inner 
product
\[\langle (a_1,\ldots ,a_n)^t,(b_1,\ldots ,b_n)^t\rangle_B =
\sum_{i=1}^n a_i^*b_i,\ \ \ \ a_i,b_i \in B.\]
Then, by \cite[Theorem 3.1]{B} and the part of the proof found on~\cite[p.\ 266]{B},
there exist nets of adjointable contractive maps
\[H \overset{\phi_{\alpha}}{\longrightarrow} 
C_{n(\alpha)}(B)\overset{\psi_{\alpha}}{\longrightarrow} H,\ \ \ \ \alpha \in 
N,\]
such that
\[\lim_{\alpha}\,\|\psi_{\alpha}(\phi_{\alpha}(h))-h\|=0,\ \ \ \ h \in H.\]
These maps induce a diagram
\begin{equation}\label{eq2.10}
\KEu(H) \overset{\widetilde\phi_{\alpha}}{\longrightarrow} 
\KEu(C_{n(\alpha)}(B))\overset{\widetilde\psi_{\alpha}}{\longrightarrow} 
\KEu(H),\ 
\ \ \ \ \alpha \in N,
\end{equation}
of complete contractions given by
\[\widetilde\phi_{\alpha}(S)=\phi_{\alpha}S\phi_{\alpha}^*,\ \ \ \ 
\widetilde\psi_{\alpha}(T)=\psi_{\alpha}T\psi_{\alpha}^*,\ \ \ \ \alpha \in 
N,\]
for $S \in \KEu(H)$ and $T \in \KEu(C_{n(\alpha)}(B))$, whose compositions 
converge to ${\mathrm{id}}_{\KEu(H)}$ in the point norm topology. The relations
\[\phi_{\alpha}\theta_{h_1,h_2}\phi_{\alpha}^*=\theta_{\phi_{\alpha}(h_1), 
\phi_{\alpha}(h_2)},\ \ \ \ \ 
\psi_{\alpha}\theta_{k_1,k_2}\psi_{\alpha}^*=\theta_{\psi_{\alpha}(k_1), 
\psi_{\alpha}(k_2)},\]
for $h_i\in H$, $k_i \in C_{n(\alpha)}(B)$, are easy to check and show that 
$\widetilde\phi_{\alpha}$ and $\widetilde\psi_{\alpha}$ have the appropriate 
ranges in (\ref{eq2.10}). 
Now $\KEu(C_{n(\alpha)}(B))$ is the matrix algebra ${\mathbb 
M}_{n(\alpha)}(B)$ over 
$B$, for which $\Lambda(\KEu(C_{n(\alpha)}(B)))=\Lambda({\mathbb 
M}_{n(\alpha)}(B))=\Lambda(B)$.

Given $S_1,\ldots ,S_r \in  \KEu(H)$ and $\varepsilon >0$, choose $\alpha \in 
N$ so large that
\begin{equation}\label{eq2.11} 
\|\widetilde\psi_{\alpha}(\widetilde\phi_{\alpha}(S_i))-S_i\|<\varepsilon /2, 
\ \ \ \ \  
1\leq i \leq r,
\end{equation}
and then choose a finite rank map $\mu \colon \KEu(C_{n(\alpha)}(B))\to 
\KEu(C_{n(\alpha)}(B))$ such that $\|\mu \|_{cb} \leq \Lambda(B)$ and
\begin{equation}\label{eq2.12}
\|\mu(\widetilde\phi_{\alpha}(S_i))-\widetilde\phi_{\alpha}(S_i)\| 
<\varepsilon /2,
\ \ \ \ \ 1\leq i \leq r.
\end{equation}
Then $\widetilde\psi_{\alpha}\mu\widetilde\phi_{\alpha}\colon 
\KEu(H) 
\to \KEu(H)$, the composition of the diagram
\[\KEu(H) \overset{\widetilde\phi_{\alpha}}{\longrightarrow} 
\KEu(C_{n(\alpha)}(B))\overset{\mu}{\longrightarrow} 
\KEu(C_{n(\alpha)}(B))\overset{\widetilde\psi_{\alpha}}{\longrightarrow} 
\KEu(H),\]
 has finite rank with completely bounded norm at 
most $\Lambda(B)$, and the triangle inequality gives
\[\|\widetilde\psi_{\alpha}\mu\widetilde\phi_{\alpha}(S_i)-S_i\| <\varepsilon,
\ \ \ \ \ 1 \leq i \leq r,\]
from (\ref{eq2.11}) and (\ref{eq2.12}).
This shows that $\Lambda(\KEu(H)) \leq \Lambda(B)$.
\end{proof}

We are now able to state and prove the main result of the paper.

\begin{thm}\label{main}
Let $B$ be a $C^*$--algebra and let $H$ be a Hilbert $B,B$--bimodule
such that $\{\langle h_1,h_2\rangle_B\colon h_1,h_2\in H\}$ generates $B$. 

Consider the extended Cuntz--Pimsner $C^*$--algebra $E(H)$.
Then $\Lambda(E(H))=\Lambda(B)$.
In particular, $E(H)$ has the completely bounded approximation property if and 
only if $B$ does.
\end{thm}
\begin{proof}
We will follow quite closely the proof of~\cite[Thm.\ 3.1]{DS}.
Since $B$ is contained as a $C^*$--subalgebra of $E(H)$ which is the image of 
a conditional
expectation $E(H)\to B$, we have $\Lambda(B)\le\Lambda(E(H))$.

Let $\Ht=H\oplus B$ and identify $H$ with the submodule $H\oplus 
0\subseteq\Ht$.
Then (see~\cite{P}), $E(H)$ is contained in $E(\Ht)$ in the obvious way.
There is a projection $P:\Ht\to H$ that commutes with the left action of $B$,
and using Proposition~\ref{prop:ST} to take tensor products of copies of $P$ 
yields a projection
$F(P):\FEu(\Ht)\to \FEu(H)$, compression with respect to which gives a 
conditional 
expectation $E(\Ht)\to E(H)$.
Therefore, $\Lambda(E(H))\le\Lambda(E(\Ht))$.

We will show that $\Lambda(E(\Ht))\le\Lambda(B)$.
By~\cite[Claim 3.4]{DS}, $E(\Ht)$ is isomorphic to the crossed product 
$A\rtimes_\Psi\Nats$ of
a certain $C^*$--subalgebra $A\subseteq E(\Ht)$ by an injective endomorphism 
$\Psi$, that is
given by $\Psi(a)=LaL^*$ for an isometry $L\in E(\Ht)$.
As described by Cuntz~\cite{C} and Stacey~\cite{S} (see also the discussion on 
p.\ 432 of~\cite{DS}),
the crossed product $A\rtimes_\Psi\Nats$ is isomorphic to a corner 
$p(\At\rtimes_\alpha\Ints)p$
of the crossed product of a $C^*$--algebra $\At$ by an automorphism $\alpha$, 
where $\At$ is the 
inductive limit $C^*$--algebra of the system
\begin{equation}\label{eq:APsi}
A\overset{\Psi}{\to}A\overset{\Psi}{\to}A\overset{\Psi}{\to}\cdots.
\end{equation}
We have 
$\Lambda(E(\Ht))=\Lambda(p(\At\rtimes_\alpha\Ints)p)\le\Lambda(\At\rtimes_\alpha\Ints)$.
From~\cite[Thm.\ 3.4]{SS}, we have 
$\Lambda(\At\rtimes_\alpha\Ints)=\Lambda(\At)$.
Letting $\sigma:A\to A$ be $\sigma(a)=L^*aL$, we have that
$\sigma$ is a completely positive left inverse of $\Psi$.
Thus $\Psi\circ\sigma$ is a conditional expectation from $A$ onto $\Psi(A)$
and Proposition~\ref{prop:Lind} applies to the inductive limit~\eqref{eq:APsi}.
We conclude that $\Lambda(\At)=\Lambda(A)$.

From the proof of~\cite[Claim 3.5]{DS}, we have an increasing chain
\begin{equation}\label{eq:AnA}
B=A_0\subseteq A_1\subseteq A_2\subseteq\cdots\subseteq A
\end{equation}
of $C^*$--subalgebras of $A$, such that $\bigcup_{n=1}^\infty A_n$ is dense in 
$A$.
Moreover, for each $n\ge1$ there is an ideal $I_n\subseteq A_n$ such that
$A_n/I_n$ is isomorphic to $A_{n-1}$; this isomorphism followed by the 
inclusion $A_{n-1}\hookrightarrow A_n$
splits the
short exact sequence
\[
0\to I_n\to A_n\to A_n/I_n\to 0.
\]
Finally, the ideal $I_n$ is isomorphic to $\KEu(H^{\otimes_Bn})$.
By Proposition~\ref{prop:KH}, $\Lambda(I_n)\le\Lambda(B)$.
Using Theorem~\ref{thm:Lses} and proceeding by induction, one shows that
$\Lambda(A_n)\le\Lambda(B)$
for all $n\ge0$.
The quotient map $A_n\to A_n/I_n$ followed by the isomorphism $A_n/I_n\to 
A_{n-1}$ is
a left inverse for the inclusion $A_{n-1}\hookrightarrow A_n$.
Thus, Proposition~\ref{prop:Lind} applies to the chain~\eqref{eq:AnA} and we 
obtain $\Lambda(A)\le\Lambda(B)$.

These inequalities, applied sequentially, lead to 
$\Lambda(E(H))\le\Lambda(B)$.
\end{proof}

\section{Concluding remarks}

An interesting open problem is whether taking reduced (amalgamated) free 
products
preserves the class of $C^*$--algebras that possess the CBAP, or, for example, 
the class
of $C^*$--algebras having Haagerup invariant equal to $1$.
It was shown in~\cite[Prop.\ 4.2]{DS} that a reduced free product 
$C^*$--algebra $A$,
\begin{equation}\label{eq:Afp}
(A,\phi)=(A_1,\phi_1)*(A_2,\phi_2),
\end{equation}
can be embedded into $E(H)$ for a particular Hilbert bimodule $H$ over
$A_1\otimes_{\min} A_2$.
It was also shown in~\cite[Prop.\ 5.1]{DS} that a reduced amalgamated free 
product $C^*$--algebra $A$,
\[
(A,\phi)=(A_1,\phi_1)*_B(A_2,\phi_2),
\]
can be realized as a quotient of a subalgebra of $E(H')$ for a particular 
Hilbert bimodule $H'$ over
$A_1\oplus A_2$.
This, together with the result from~\cite{DS} that $E(H)$ is exact when $H$ is 
a Hilbert bimodule
over an exact $C^*$--algebra, yielded a new proof that the class of exact 
$C^*$--algebras is closed
under taking reduced (amalgamated) free products.
(See~\cite{D} for the first proof of this fact.)

However, this paper's Theorem~\ref{main} does not answer in a similar way  the 
question about the CBAP
for free products, (at least not obviously), because the CBAP does not 
automatically pass to subalgebras 
or to quotients.
Moreover, the copy of the free product $C^*$--algebra $A$ from~\eqref{eq:Afp} 
in $E(H)$ exhibited
in~\cite{DS} is not in general the image of a conditional expectation $E(H)\to 
A$.
To see this, first note from~\cite{G} that
$E(H)$ is a nuclear $C^*$--algebra
whenever $H$ 
is a Hilbert bimodule over a nuclear $C^*$--algebra.
(This also follows readily from the proof of exactness found in~\cite{DS}, combined with the conditional expectation
$E(\Ht)\to E(H)$ used in the proof of Theorem~\ref{main}.)

There are many known examples when $A_1$ and $A_2$ in~\eqref{eq:Afp} are 
nuclear, but their free product $A$
is not nuclear.
In these cases, $E(H)$ is nuclear and $A$, therefore, cannot be the image of a 
conditional expectation $E(H)\to A$.

\bibliographystyle{plain}

\end{document}